\documentclass[12pt]{article}
\usepackage{amsfonts}
\usepackage{amssymb}
\usepackage{latexsym}
\usepackage{amsmath}
\usepackage{bibmods}
\usepackage{graphicx}
\usepackage{color}
\usepackage[arrow,matrix,curve]{xy}
\usepackage{amssymb}
\usepackage{amscd}
\usepackage{color}

\newtheorem{lemma}{\indent Lemma}[section]
\newtheorem{theorem}[lemma]{\indent Theorem}

\newtheorem{remark}[lemma]{\indent Remark}
\newcommand{\cx}{{}{\scriptstyle{\mathcal X}}}
\newcommand{\wt}{\widetilde}
\begin{document}

\begin{center}
{\Large\textbf{ON POINCAR\'E, FRIEDRICHS AND KORNS INEQUALITIES ON DOMAINS AND HYPERSURFACES}\footnote{This work was supported by the Shota Rustaveli Georgian National Science Foundation No. GNSF/DI/10/5-101/12, Contract No. 13/14.}}
\vspace{5mm}
 \end{center}
\begin{center}

\textbf{R. Duduchava}

\vspace{2mm}

I. Javakhishvili Tbilisi State University, Andrea Razmadze
Mathematical Institute,  Tamarashvili str. 6, Tbilisi 0177, Georgia;
\texttt{roland.duduchava@tsu.ge}
 \end{center}

\thispagestyle{empty}

\begin{abstract}
The celebrated Poincar\'e and Friedrichs inequalities estimate the $\mathbb{L}_p$-norm of a function by the $\mathbb{L}_p$-norm of the gradient. We prove the Poincar´e inequality for a domain $\Omega\subset \mathbb{R}^n$ and for a hypersurface $\mathcal{C}\subset\mathbb{R}^n$ based on open mapping theorem of Banach only. For a cylinder which has a hypersurface as a base, is prove stronger inequality, involving only the surface derivatives. Similar inequalities for the uniform $C$-norm are proved as well. We also estimate $\mathbb{H}^m_p$-norm of functions  prove inequalities for some generalizations of the mentioned inequalities.

We also prove Poincar\'e-Korns and Friedrichs-Korns inequalities for vector-func\-ti\-ons estimating the $\mathbb{L}_p$-norm of a function by the $\mathbb{L}_p$-norm of the deformation tensor only on domains and on hypersurfaces. The proofs are based on the paper \cite{Du10} of the author on Korns inequalities. And again, the norm of the function in a cylinder is estimated by is the deformation tensor on the base of the cylinder.
\end{abstract}

 %
 %
 \section*{Introduction}
 \setcounter{equation}{0}
 \label{sec0}

Let $1\leqslant p \leqslant\infty$ and $\Omega$ be a bounded connected open subset of the $n$-dimensional Euclidean space $\mathbb{R}^n$ with a Lipschitz boundary (a domain with the {\em uniform cone property}). Then there exists a constant C, depending only on $\Omega$ and $p$ such that for every function $\varphi$ in the Sobolev space $\mathbb{W}^1_p(\Omega)$ the celebrated Poincar\'e inequality holds
 \begin{eqnarray}\label{e0.1P}
\|\varphi - \varphi_\Omega\big|\mathbb{L}_p(\Omega)\|\leqslant C\|\nabla\varphi\big|
     \mathbb{L}_p(\Omega)\|,
 \end{eqnarray}
where
 \begin{eqnarray}\label{e0.2}
\varphi_\Omega:=\frac1{{\rm mes}\,\Omega}\int_\Omega\varphi(y)dy
 \end{eqnarray}
is the average value of $\varphi$ over $\Omega$. Here ${\rm mes}\,\Omega$ stands for the Lebesgue measure of the domain $\Omega$ and the constant $C$ depends on $\Omega$ and $p$ only. When $\Omega$ is a ball, the above inequality is called a {\em Poincar\'e inequality}, while for more general domains $\Omega$ inequality \eqref{e0.1P} is known as a {\em Sobolev inequality} (cf., e.g., \cite{DL90}).

Let $\mathcal{M}_0$ be a subset of the closed domain $\mathcal{M}_0\subset\overline\Omega$ of co-dimension 1 and have non-trivial measure ${\rm mes}\,\mathcal{M}_0\not=0$ (can be a non-trivial part of the boundary). Let $\varphi^+$ denote the trace of $\varphi$ on $\mathcal{M}_0$. The following
 \begin{eqnarray} \label{e0.1F}
\hskip10mm\|\varphi\big|\mathbb{L}_p(\Omega)\|\leqslant C\left[\|\nabla\varphi\big|
     \mathbb{L}_p(\Omega)\|^p+\left|\int_{\mathcal{M}_0}\varphi^+(\cx)\,
     d\sigma\right|^p\right]^{1/p},\qquad \varphi\in\mathbb{W}^1_p(\Omega)
 \end{eqnarray}
is known as {\em Friedrichs inequality} for $\mathcal{M}_0=\partial\Omega$, $p=2$ (see \cite[Theorem 6.28.2]{Tr72}, \cite[Theorem 4.1.7]{HW08}).

If $\mathcal{M}_0$ is the same as in \eqref{e0.1F}, the next inequality
 \begin{eqnarray}\label{e0.3}
\|\varphi\big|\mathbb{L}_p(\Omega)\|\leqslant C\|\nabla\varphi\big|\mathbb{L}_p(\Omega)\|
 \end{eqnarray}
for a function $\varphi\in\widetilde{\mathbb{W}}^1_p(\Omega,\mathcal{M}_0)$ which vanish on $\mathcal{M}_0$,
is a variant of inequalities \eqref{e0.1P}, \eqref{e0.1F} (see \cite[Theorem 6.28.2]{Tr72}, \cite[Theorem 4.1.7]{HW08} and \cite[Theorem 7.6, Theorem 7.7]{Wl87}).

The inequalities \eqref{e0.1F} and \eqref{e0.3} hold, of course, if $\mathcal{M}_0$ is a subdomain of $\Omega$.

In contrast to \eqref{e0.1P}, in inequalities \eqref{e0.3} and \eqref{e0.1F} the domain $\Omega$ can also be unbounded (might have an infinite measure), provided ${\rm mes}\,\mathcal{M}_0<\infty$ in \eqref{e0.1F}.

Moreover, for a cylinder $\Omega:=\mathcal{C}\times[a,b]$ with a base $\mathcal{C}$ which is a hypersurface in $\mathbb{R}^n$, we prove a stronger inequality, namely the following
 \begin{eqnarray}\label{e0.5}
&\hskip5mm\|\varphi\big|\mathbb{L}_p(\Omega)\|\leqslant C\left[\|\nabla_\mathcal{C}\varphi
     \big|\mathbb{L}_p(\Omega)\|^p+\left|\displaystyle\int_{\mathcal{M}_0}\varphi^+(\cx)\,d\sigma
     \right|^p\right]^{1/p},\qquad\varphi\in\mathbb{W}^1_p(\mathcal{C}),\\[2mm]
\label{e0.5a}
&\|\varphi\big|\mathbb{L}_p(\Omega)\|\leqslant C\|\nabla_\mathcal{C}\varphi\big|
     \mathbb{L}_p(\Omega)\|,\qquad \varphi\in\widetilde{\mathbb{W}}^1_p (\Omega,\mathcal{M}_0),
 \end{eqnarray}
where $\nabla_\mathcal{C}=(\mathcal{D}_1,\ldots,\mathcal{D}_n)^\top$ is the surface gradient and $\mathcal{D}_1,\ldots,\mathcal{D}_n$ are the Gunter's derivatives (see \S\, 1),  and $\varphi\in\widetilde{\mathbb{W}}^1_p (\Omega,\mathcal{M}_0)$ vanishes on a $(n-1)$-dimensional strip $\mathcal{M}_0:=\Gamma_0\times[a,b]$ with $\Gamma_0\subset\overline{\mathcal{C}}$-a $(n-2)$-dimensional subset of $\overline{\mathcal{C}}$ (can be a piece of the boundary $\partial\mathcal{C}$). The inequality \eqref{e0.5} is remarkable, because contains only the surface derivatives and does not contains the derivative with respect to the variable $t\in[a,b]$ transversal to the surface $\mathcal{C}$.

For a cylinder $I_\omega:=\omega\times I$, $I:=[a,b]$, with a flat base $\omega\subset\mathbb{R}^{n-1}$, the inequalities \eqref{e0.5} and  \eqref{e0.5a} have the form
 \begin{eqnarray}\label{e0.5b}
&&\hskip-10mm\|\varphi\big|\mathbb{L}_p(I_\omega)\|\leqslant C\left[\|\nabla_\omega\varphi
     \big|\mathbb{L}_p(I_\omega)\|^p+\|\varphi\big|\mathbb{L}_p(\Gamma_0\times I)\|^p\right]^{1/p},\\
&&\hskip60mm\varphi\in\mathbb{W}^1_p(I_\omega),\nonumber\\[2mm]
\label{e0.5c}
&&\hskip-10mm\|\varphi\big|\mathbb{L}_p(I_\omega)\|\leqslant C\|\nabla_\omega
     \varphi\big|\mathbb{L}_p(I_\omega)\|,\qquad \varphi\in\widetilde{\mathbb{W}}^1_p (I_\omega,\Gamma_0\times I),
 \end{eqnarray}
where $\nabla_\omega(\mathbf{U})$ is the gradient in $\omega$ (in $\mathbb{R}^{n-1}$) and  contains only $(n-1)$ derivatives.

Poincar\'e and Friedrichs inequalities also hold for a smooth surfaces
 \begin{eqnarray}\label{e0.6}
&\|\varphi - \varphi_\mathcal{C}\big|\mathbb{L}_p(\mathcal{C})\|\leqslant C\|\nabla_\mathcal{C}\varphi\big|
     \mathbb{L}_p(\mathcal{C})\|,\qquad \varphi\in\mathbb{W}^1_p(\mathcal{C}),\\[2mm]
\label{e0.6F}
&\hskip5mm\|\varphi\big|\mathbb{L}_p(\mathcal{C})\|\leqslant C\left[\|\nabla_\mathcal{C}\varphi\big|
     \mathbb{L}_p(\mathcal{C})\|^p+\left|\displaystyle\int_{\Gamma_0}\varphi^+(\cx)\,d\sigma\right|^p
     \right]^{1/p},\quad \varphi\in\mathbb{W}^1_p(\mathcal{C})\\[2mm]
\label{e0.7}
&\|\varphi\big|\mathbb{L}_p(\mathcal{C})\|\leqslant C\|\nabla_\mathcal{C}\varphi\big|
     \mathbb{L}_p(\mathcal{C})\|,\qquad \varphi\in\widetilde{\mathbb{W}}^1_p(\mathcal{C},\Gamma_0),
 \end{eqnarray}
where $\varphi_\mathcal{C}$ denotes the average value of $\varphi$ over $\mathcal{C}$:
 \begin{eqnarray}\label{e0.4}
\varphi_\mathcal{C}:=\frac1{{\rm mes}\,\mathcal{C}}\int_\mathcal{C}\varphi(y)d\sigma.
 \end{eqnarray}
$\Gamma_0$ is a subset of the closed surface $\Gamma_0\subset\overline{\mathcal{C}}$ of co-dimension 1 and has non-trivial measure ${\rm mes}\,\Gamma_0\not=0$ ($\Gamma_0$ can be a non-trivial part of the boundary).

The inequalities \eqref{e0.6} and \eqref{e0.7} hold, of course, if $\Gamma_0$ is a subsurface of $\mathcal{C}$.

The inequality \eqref{e0.6} holds for surfaces of finite measure, while the inequality \eqref{e0.7} does not needs such constraint and the surface $\mathcal{C}$ might have infinite measure.

The following
 \begin{eqnarray} \label{e0.14P}
&&\|\varphi\big|\mathbb{W}^\ell_p(\Omega)\|\leqslant\|\varphi\big|\mathbb{W}^m_p(\Omega)\|
     \leqslant C\left[\sum_{|\alpha|=m}\|\partial^\alpha\varphi\big|\mathbb{L}_p(\Omega)\|^p \right.\nonumber\\
&&\hskip10mm\left.+\sum_{|\beta|<m}\left|\int_{\mathcal{M}_0}(\partial^\beta\varphi)^+(\cx)\,d\sigma
     \right|^p\right]^{1/p},\qquad\varphi\in\mathbb{W}^m_p(\Omega),\\
\label{e0.14P1}
&&\|\varphi\big|\mathbb{W}^\ell_p(\mathcal{M})\|\leqslant\|\varphi\big|\mathbb{W}^m_p(\mathcal{M})\|
     \leqslant C\left[\sum_{|\alpha|=m}\|\mathcal{D}^\alpha\varphi\big|\mathbb{L}_p(\mathcal{M})\|^p \right.\nonumber\\
&&\hskip10mm\left.+\sum_{|\beta|<m}\left|\int_{\Gamma_0}(\mathcal{D}^\beta\varphi)^+(\cx)\,d\sigma
     \right|^p\right]^{1/p},\qquad\varphi\in\mathbb{W}^m_p(\mathcal{M}),\\
\label{e0.14P2}
&&\hskip-15mm\|\varphi\big|\mathbb{W}^\ell_p(\Omega)\|\leqslant\|\varphi\big|\mathbb{W}^m_p
     (\Omega)\|\leqslant C\sum_{|\alpha|=m}\|\partial^\alpha\varphi\big|\mathbb{L}_p(\Omega)\|^p
     \quad \varphi\in\widetilde{\mathbb{W}}^m_p(\Omega,\mathcal{M}_0),
 \end{eqnarray}
 \begin{eqnarray} 
\label{e0.14P3}
&&\hskip-15mm\|\varphi\big|\mathbb{W}^\ell_p(\mathcal{M})\|\leqslant\|\varphi\big|\mathbb{W}^m_p(\mathcal{M})\|
     \leqslant C\sum_{|\alpha|=m}\|\mathcal{D}^\alpha\varphi\big|\mathbb{L}_p(\mathcal{M})\| \quad
     \varphi\in\widetilde{\mathbb{W}}^m_p(\mathcal{M},\Gamma_0)
 \end{eqnarray}
for $\ell<m$, $m=2,3,\ldots$, generalize {\em Poincar\'e inequalities} \eqref{e0.1P} and \eqref{e0.6}, while the inequalities
 \begin{eqnarray} \label{e0.14F}
&&\hskip-15mm\|\varphi\big|\mathbb{W}^\ell_p(\Omega)\|\leqslant\|\varphi\big|\mathbb{W}^m_p(\Omega)\|
     \leqslant C\left[\sum_{|\alpha|=m}\|\partial^\alpha\varphi\big|\mathbb{L}_p(\Omega)\|^p
     +\int_{\mathcal{M}_0}|\varphi^+(\cx)|^pd\sigma\right]^{1/p}\hskip-3mm,\\
&&\hskip100mm\varphi\in\mathbb{W}^m_p(\Omega),\nonumber\\
\label{e0.14F1}
&&\hskip-15mm\|\varphi\big|\mathbb{W}^\ell_p(\mathcal{M})\|\leqslant\|\varphi\big|\mathbb{W}^m_p(\mathcal{M})\|
     \leqslant C\left[\sum_{|\alpha|=m}\|\partial^\alpha\varphi\big|\mathbb{L}_p(\mathcal{M})\|^p +\hskip-2mm\int_{\Gamma_0}\hskip-3mm|\varphi^+(\cx)|^pd\sigma\right]^{1/p}\hskip-3mm,\\
&&\hskip100mm\varphi\in\mathbb{W}^m_p(\mathcal{M})\nonumber
 \end{eqnarray}
for $\ell<m$, $m=2,3,\ldots$, generalize {\em Friedrichs inequalities} \eqref{e0.1F} and \eqref{e0.6F} (see \cite[Theorem 6.28.2]{Tr72}, \cite[Theorem 4.1.7]{HW08}).

All above inequalities hold also for the space of 1-smooth functions $C^1(\Omega)$-just replace the $\mathbb{L}_p$-norm by $\|\varphi\big|C(\mathcal{C})\|:=\max_{x\in\Omega}|\varphi(x)|$ and $\mathbb{W}^1_p$-norm by $\|\varphi\big|C^1(\mathcal{C})\|:=\max_{x\in\Omega}|\varphi(x)| +\max_{x\in \Omega}|\nabla\varphi(x)|$.. For example, the inequality \eqref{e0.1P} acquires the form
 \begin{eqnarray}\label{e0.8}
\max_{x\in\Omega}|\varphi(x) - \varphi_\Omega(x)|\leqslant C\max_{x\in\Omega}|
     \nabla\varphi(x)|.
 \end{eqnarray}
There is only one essential difference: in analogues of inequalities \eqref{e0.3}, \eqref{e0.4}, \eqref{e0.6}, \eqref{e0.14P} and \eqref{e0.14F} the sets $\mathcal{M}_0$ and $\Gamma_0$ can be one point sets.

It turned out, that for vector-functions $\mathbf{U}(x)=(U_1(x),\ldots,U_n(x))^\top$ on a domain $\Omega\subset\mathbb{R}^n$ even gradient is superfluous in the inequalities \eqref{e0.1P}, \eqref{e0.3}  and it suffices to take the deformation tensor:
 \begin{eqnarray}
\label{e0.11}
&\hskip7mm\|\mathbf{U}\big|\mathbb{L}_p(\Omega)\|\leqslant C\left[\|{\rm\bf Def}\,\mathbf{U}\big|
    \mathbb{L}_p(\Omega)\|^p+\|\mathbf{U}^+\big|\mathbb{L}_p(\mathcal{M}_0)\|^p\right]^{1/p},\quad \mathbf{U}\in\mathbb{W}^1_p(\Omega),\\[2mm]
\label{e0.12}
&\hskip7mm\|\mathbf{U}\big|\mathbb{L}_p(\mathcal{C})\|\leqslant C\left[\|{\rm\bf Def}_\mathcal{C}\mathbf{U}\big|
    \mathbb{L}_p(\mathcal{C})\|^p+\|\mathbf{U}^+\big|\mathbb{L}_p(\Gamma_0)\|^p\right]^{1/p} \quad \mathbf{U}\in\mathbb{W}^1_p(\mathcal{C}),\\[2mm]
 \label{e0.9}
&\|\mathbf{U}\big|\mathbb{L}_p(\Omega)\|\leqslant C\|{\rm\bf Def}\,\mathbf{U}\big|
    \mathbb{L}_p(\Omega)\|,\qquad \mathbf{U}\in\widetilde{\mathbb{W}}^1_p(\Omega,\mathcal{M}_0),\\[2mm]
\label{e0.10}
&\|\mathbf{U}\big|\mathbb{L}_p(\mathcal{C})\|\leqslant C\|{\rm\bf Def}_\mathcal{C}\mathbf{U}\big|\mathbb{L}_p(\mathcal{C})\|
    \qquad \mathbf{U}\in\widetilde{\mathbb{W}}^1_p(\mathcal{C},\Gamma_0),
 \end{eqnarray}
where $\mathcal{M}_0$ and $\Gamma_0$ are the same as in \eqref{e0.1F} and \eqref{e0.6F}, respectively.  ${\rm\bf Def}(\mathbf{U})$ and ${\rm\bf Def}_\mathcal{C}(\mathbf{U})$ are the domain and the surface deformation tensors, respectively (see \eqref{e1.10} and \eqref{e1.11}), and only $\displaystyle\frac{n(n+1)}2<n^2$ different linear combinations of the $n^2$ derivatives $\partial_jU_k$ (of derivatives $\mathcal{D}_jU_k$, respectively; $j,k=1,\ldots,n$) are involved.

For a cylinder $\Omega:=\mathcal{C}\times[a,b]$ with a base $\mathcal{C}$ which is a hypersurface in $\mathbb{R}^n$ and a vector-function $\mathbf{U}=(U_1,\ldots,U_n)^\top$, we prove a stronger inequality, namely the following
 \begin{eqnarray}\label{e0.15}
&\;\;\;\;\;\|\mathbf{U}\big|\mathbb{L}_p(\Omega)\|\leqslant C\left[\|{\rm\bf Def}_\mathcal{C}(\mathbf{U})
     \big|\mathbb{L}_p(\Omega)\|^p+\|\mathbf{U}\big|\mathbb{L}_p(\mathcal{M}_0)\|^p\right]^{1/p},
     \quad\mathbf{U}\in\mathbb{W}^1_p(\Omega),\\[2mm]
\label{e0.16}
&\;\;\;\;\;\;\;\|\mathbf{U}\big|\mathbb{L}_p(\Omega)\|\leqslant C\|{\rm\bf Def}_\mathcal{C}(\mathbf{U})\big|
     \mathbb{L}_p(\Omega)\|,\qquad \mathbf{U}\in\widetilde{\mathbb{W}}^1_p (\Omega,\mathcal{M}_0),
 \end{eqnarray}
where $\mathcal{M}_0:=\Gamma_0\times[a,b]$ is a strip, $\Gamma_0\subset\overline{\mathcal{C}}$. The inequality \eqref{e0.5} is remarkable, because estimates the vector-function $\mathbf{U}$, instead of  $n(n+1)$ derivatives $\mathcal{D}_jU_k$, $j=1,\ldots,n+1,\ k=1,\ldots,n$ including the transversal derivatives $\mathcal{D}_{n+1}U_k$, $k=1,\ldots,n$, by only surface deformation tensor ${\rm\bf Def}_cC(\mathbf{U})$.

For a cylinder $I_\omega:=\omega\times I$, $I:=[a,b]$, with a flat base $\omega\subset\mathbb{R}^{n-1}$, the inequalities \eqref{e0.15} and  \eqref{e0.16} have the form
 \begin{eqnarray}\label{e0.17}
&&\hskip-10mm\|\mathbf{U}\big|\mathbb{L}_p(I_\omega)\|\leqslant C\left[\|{\rm\bf Def}_\omega(\mathbf{U})
     \big|\mathbb{L}_p(I_\omega)\|^p+\|\mathbf{U}\big|\mathbb{L}_p(\Gamma_0\times I)\|^p\right]^{1/p},\\
&&\hskip65mm\mathbf{U}\in\mathbb{W}^1_p(I_\omega),\nonumber\\[2mm]
\label{e0.18}
&&\hskip-10mm\|\mathbf{U}\big|\mathbb{L}_p(I_\omega)\|\leqslant C\|{\rm\bf Def}_\omega(\mathbf{U})\big|\mathbb{L}_p
     (I_\omega)\|,\qquad \mathbf{U}\in\widetilde{\mathbb{W}}^1_p (I_\omega,\Gamma_0\times I),
 \end{eqnarray}
where ${\rm\bf Def}_\omega(\mathbf{U})$ is the deformation tensor in $\omega$ (in $\mathbb{R}^{n-1}$) and  contains only $\displaystyle\frac{n(n-1)}2$ derivatives.

The inequalities \eqref{e0.11}-\eqref{e0.10} follow from Korns inequalities and we call:  \eqref{e0.11}-\eqref{e0.12} {\em Friedrichs-Korns inequalities} and  \eqref{e0.9}-\eqref{e0.10}  {\em Poincar\'e-Korns inequalities}.

 %
 %
\section{Auxiliaries}
 \label{sec1}
\setcounter{equation}{0}

Throughout the present paper we will assume that $\mathcal{C}$ be a sufficiently smooth hypersurface in $\mathbb{R}^n$ with the Lipschitz boundary $\Gamma:=\partial\mathcal{C}$  (a surface with the {\em uniform cone property}), defined by a real valued smooth function
 \begin{equation}\label{e1.1}
\mathcal{C}=\Big\{\cx\in \Omega\;:\; \Psi_\mathcal{C}(\cx)=0\Big\},
 \end{equation}
which is regular $\nabla\,\Psi_\mathcal{C}(\cx)\not=0$. The normalized  gradient
 \begin{equation}\label{e1.2}
\boldsymbol{\nu}(\cx):=\displaystyle\frac{\nabla\Psi_\mathcal{C}(\cx)}{\big|\nabla\Psi_\mathcal{C}(\cx)\big|}\, ,\qquad \cx\in\mathcal{C}
 \end{equation}
defines the {\em unit normal vector field} on $\mathcal{C}$.

The collection of the tangential {\em G\"unter's derivatives} are defined as follows (cf. \cite{Gu53,KGBB79,DMM06,Du10,Du11})
 \begin{eqnarray}\label{e1.3}
\qquad\qquad\mathcal{D}_j:=\partial_j-\nu_j(\cx)\partial_{\boldsymbol{\nu}}=\partial_{\mathbf{d}\,^j},\qquad
     \nu_j(\cx):=\displaystyle\frac{\partial_j\Psi_\mathcal{C}(\cx)}{\big|\nabla\Psi_\mathcal{C}(\cx)\big|}\,\quad
     j=1,\ldots,n,
 \end{eqnarray}
where
 \begin{equation}\label{e1.4}
 \mathbf{e}^1=(1,0,\ldots,0)^\top\, ,\ldots, \mathbf{e}^n=(0,\ldots,0,1)^\top
 \end{equation}
is the natural basis in $\mathbb{R}^n$ and $\partial_{\boldsymbol{\nu}}:=\sum_{j=1}^n\nu_j\partial_j$ denotes the normal derivative. For each $1\leq j\leq n$, the first-order differential operator ${\mathcal{D}}_j=\partial_{ \mathbf{d}\,^j}$ is the directional derivative along the tangential vector
 \begin{equation}\label{e1.5}
\mathbf{d}\,^j:=\pi_\mathcal{S}\,\mathbf{e}^j, \quad \langle\boldsymbol{\nu}(x),\mathbf{d}\,^j(x)\rangle\equiv0,\quad
     \sum_{j=1}^n\nu_k\mathbf{d}\,^k=0, \quad j=1,\ldots,n,
 \end{equation}
the projection of $\mathbf{e}^j$ on the space of tangential vector fields to $\mathcal{S}$.

The {\em surface gradient} $\nabla_\mathcal{S} \varphi$ is the collection of the G\"unter's derivatives
 \begin{eqnarray}\label{e1.6}
\nabla_\mathcal{S}\varphi:=(\mathcal{D}_1\varphi,\ldots,\mathcal{D}_n\varphi)^\top
 \end{eqnarray}
and is an equivalent form of the surface gradient defined in the differential geometry by means of covariant metric tensor (see \cite{DMM06,Du10,Du11}). The next Lemma \ref{l1.1} was proved in \cite[Lemma 1.2]{Du10}.
%
\begin{lemma}\label{l1.1}
For $\varphi\in C^1(\mathcal{S})$ the surface gradient vanishes $\nabla_\mathcal{S}\varphi\equiv0$ if and only if $\varphi(\cx)\equiv{\rm const}$.
\end{lemma}

$\mathbb{W}^1_p(\Omega)$ and $\mathbb{W}^1_p(\mathcal{C})$, $1<p<\infty$, denote the Sobolev spaces on a domain $\Omega\subset\mathbb{R}^n$ and the surface $\mathcal{C}$ endowed with the norm:
 \begin{equation}\label{e1.7}
\|\varphi\,\big|\,\mathbb{W}_p^1(\Omega)\,\|:=\left[\|\varphi\,\big|\,\mathbb{L}_p(\Omega)\,\|
    +\sum\limits_{j=1}^n\|\partial_j\varphi\,\big|\,\mathbb{L}_p(\Omega)\|_p\right]^{1/p}
 \end{equation}
and, respectively,
 \begin{equation}\label{e1.8}
\|\varphi\,\big|\,\mathbb{W}_p^m(\mathcal{C})\,\|:=\left[\|\varphi\,\big|\,\mathbb{L}_p(\Omega)\,\|
    +\sum\limits_{j=1}^n\|\mathcal{D}_j\varphi\,\big|\,\mathbb{L}_p(\mathcal{C})\|_p\right]^{1/p}.
 \end{equation}

Let us define the space $\wt{\mathbb{W}}^1_p(\Omega,\mathcal{M}_0)$ for a domain $\Omega\subset\mathbb{R}^n$ with a Lipshitz boundary $\mathcal{M}:=\partial\Omega$ and a subsurface $\mathcal{M}_0\subset\mathcal{M}$-of non-zero measure as the closure in $\mathbb{W}^1_p(\Omega)$ of the set $C^\infty(\Omega,\mathcal{M}_0)$ of smooth functions $\varphi(x)$ which have vanishing traces on $\mathcal{M}_0$, i.e. $\varphi^+(\cx)=0$ for all $\cx\in\mathcal{M}_0$.
The space $\wt{\mathbb{W}}^1_p(\Omega,\mathcal{M}_0)$ inherits the standard norm $\|\varphi\,\big|\,\mathbb{W}_p^1(\Omega)\,\|$ from the space $\mathbb{W}^1_p(\Omega)$ (see \eqref{e1.7}).

If $\mathcal{C}$ is a subsurface of a closed surface $\mathcal{S}$ without boundary, $\wt{\mathbb{W}}^1_p(\mathcal{C})$ denotes the space of functions $\varphi\in\mathbb{W}^1_p(\mathcal{S})$, supported in $\overline{\mathcal{C}}$. Let $\mathcal{C}^c=\mathcal{S}\setminus\overline{\mathcal{C}}$ be the complemented surface with the common boundary $\partial\mathcal{C}=\partial\mathcal{C}^c=\Gamma$; The notation $\mathbb{W}^1_p(\mathcal{C})$ is used for the factor space $\mathbb{W}^1_p(\mathcal{C})/\wt{\mathbb{W}}^1_p(\mathcal{C}^c)$. The space $\mathbb{W}^s_p(\mathcal{C})$ can also be interpreted as the space of restrictions $r_\mathcal{C}\varphi:=\varphi\big|_\mathcal{C}$ of all functions  $\varphi\in\mathbb{W}^1_p(\mathcal{S})$ to the subsurface $\mathcal{C}$.

Similarly are defined the spaces $\wt{\mathbb{W}}^1_p(\Omega)$ and $\mathbb{W}^1_p(\Omega)$ for a domain $\Omega\subset\mathbb{R}^n$.

We refer to \cite{Tr95,Du10,Du11} for details about these spaces.

For an n-vector-function $\mathbf{U}(x)=(U_1(x),\ldots,U_n(x))^\top$ on a domain in the Euclidean space $\Omega\subset\mathbb{R}^n$ the deformation tensor reads
 \begin{eqnarray}\label{e1.10}
{\rm\bf Def}(\mathbf{U})=\big[\mathbf{D}_{jk}(\mathbf{U})\big]_{n\times n}\, ,\qquad
   \mathbf{D}_{jk}(\mathbf{U}):=\frac12\left[\partial_jU_k+\partial_kU_j\right].
 \end{eqnarray}
 The following form of the important deformation (strain) tensor on a surface $\mathcal{C}$ was identified in \cite{DMM06}:
 \begin{eqnarray}\label{e1.11}
&\;\;\;\;\;\;\;\;{\rm\bf Def}_\mathcal{C}(\mathbf{U})
    =\big[\mathfrak{D}_{jk}(\mathbf{U})
    \big]_{n\times n}\, ,\quad\mathbf{U}=\sum_{j=1}^nU_j\mathbf{d}^j\in\mathcal{V}( \mathcal{C}), \quad j,k=1,\ldots,n, \\
&\mathfrak{D}_{jk}(\mathbf{U}):=\frac12\bigl[(\mathcal{D}^\mathcal{S}_j
    \mathbf{U})_k+(\mathcal{D}^\mathcal{S}_k\mathbf{U})_j\bigr]=\frac12\Bigl[{
    \mathcal{D}}_kU_j+{\mathcal{D}}_jU_k+\sum\limits_{m=1}^nU_m\mathcal{D}_m
    \big(\nu_j\nu_k\big)\Bigr]\, , \nonumber
 \end{eqnarray}
where $(\mathcal{D}^\mathcal{S}_j\mathbf{U})_k:=\langle\mathcal{D}^\mathcal{S}_j\mathbf{U},\mathbf{e}^k\rangle$ and $\mathcal{V}(\mathcal{C})$ is the linear space of all tangential vectors-functions to the surface $\mathcal{C}$.

A vector $\mathbf{U}\in\mathbb{W}^1_p(\Omega)$ is called a {\em rigid motion} if ${\rm\bf Def}(\mathbf{U})=0$ and a vector $\mathbf{V}\in\mathbb{W}^1_p(\mathcal{C})$ is called a {\em Killings vector field} on the surface $\mathcal{C}$ if ${\rm\bf Def}_\mathcal{C}(\mathbf{V})=0$.

The next Theorem \ref{t1.2} ({\em Korns I inequality for domains ``without boundary con-\linebreak\-di\-tion"}) is well known (see \cite{Ci00} for a simple proof when $p=2$, $m=1$ and see \cite[Theorem 2.3]{Du10} for a general case).

Theorem \ref{t1.3} is proved in \cite[Theorem 2.3]{Du10}. P. Ciarlet proved it in \cite{Ci00} for the case $p=2$, $m=1$, manifold without boundary, for curvilinear coordinates and covariant derivatives.
 %
 \begin{theorem}\label{t1.2} Let $1<p<\infty$, $\Omega\subset\mathbb{R}^n$ be a domain with the
Lip\-shitz boundary and
 \begin{equation}\label{e1.12}
\big\|{\rm\bf Def}(\mathbf{U})\big|\mathbb{L}_p(\Omega)\big\|:=\left[\sum_{j,k=1}^n
    \big\|\mathbf{D}_{jk}(\mathbf{U})\big|\mathbb{L}_p(\Omega)\big\|^p\right]^{1/p},\qquad \mathbf{U}\in\mathbb{W}^1_p(\Omega).
 \end{equation}
Then the inequality
 \begin{eqnarray}\label{e1.13}
\big\|\mathbf{U}\big|\mathbb{W}^1_p(\Omega)\big\|\leq M\left[\big\|\mathbf{U}\big|
     \mathbb{L}_p(\Omega)\big\|^p+\big\|{\rm\bf Def}(\mathbf{U})\big|\mathbb{L}_p(\Omega)
     \big\|^p\right]^{1/p}
 \end{eqnarray}
holds with some constant $M>0$ or, equivalently, the equality
 \begin{eqnarray}\label{e1.14}
\big\|\mathbf{U}\big|\mathbb{W}^1_p(\Omega)\big\|_0:=\left[\big\|\mathbf{U}\big|
     \mathbb{L}_p(\Omega)\big\|^p+\big\|{\rm\bf Def}(\mathbf{U})\big|\mathbb{L}_p(\Omega)\big\|^p\right]^{1/p}
 \end{eqnarray}
defines an equivalent norm on the space $\mathbb{W}^1_p(\Omega)$.

A rigid motion $\mathbf{U}$, ${\rm\bf Def}(\mathbf{U})=0$, has the unique continuation property: if
$\mathbf{U}(x)=0$ on a set $\mathcal{M}_0$ described in \eqref{e0.1F}, than $\mathbf{U}(x)=0$ everywhere on $\Omega$.
 \end{theorem}
 %
 \begin{theorem}\label{t1.3} Let $1<p<\infty$, $\mathcal{C}\subset\mathbb{R}^n$ be a
 Lip\-shitz hypersurface with or without boundary and (see \eqref{e1.11} for the deformation tensor ${\rm\bf Def}_\mathcal{C}(\mathbf{V})$)
 \begin{equation}\label{e1.15}
\big\|{\rm\bf Def}_\mathcal{C}(\mathbf{V})\big|\mathbb{L}_p(\mathcal{C})\big\|:=\left[\sum_{j,k=1}^n
    \big\|\mathfrak{D}_{jk}(\mathbf{V})\big|\mathbb{L}_p(\mathcal{C})\big\|^p\right]^{1/p},
    \qquad \mathbf{V}\in\mathbb{W}^1_p(\mathcal{C}).
 \end{equation}
Then the inequality
 \begin{eqnarray}\label{e1.16}
\hskip-5mm\big\|\mathbf{V}\big|\mathbb{H}^1_p(\mathcal{C})\big\|\leq M\left[\big\|\mathbf{V}\big|
     \mathbb{L}_p(\mathcal{C})\big\|^p+\big\|{\rm\bf Def}_\mathcal{C}(\mathbf{V})\big|\mathbb{L}_p(\mathcal{C})
     \big\|^p\right]^{1/p}
 \end{eqnarray}
holds with some constant $M>0$ or, equivalently, the equality
 \begin{eqnarray}\label{e1.17}
\hskip-5mm\big\|\mathbf{V}\big|\mathbb{W}^1_p(\mathcal{C})\big\|_0:=\left[\big\|\mathbf{V}\big|
     \mathbb{L}_p(\mathcal{C})\big\|^p+\big\|{\rm\bf Def}_\mathcal{C}(\mathbf{V})\big|\mathbb{L}_p(\mathcal{C})
     \big\|^p\right]^{1/p}
 \end{eqnarray}
defines an equivalent norm on the space $\mathbb{W}^1_p(\mathcal{S})$.

A Killings vector field $\mathbf{V}$, ${\rm\bf Def}_\mathcal{C}(\mathbf{V})=0$, has the unique continuation property: if $\mathbf{V}(x)=0$ on a set $\Gamma_0$ described in \eqref{e0.6F}, than $\mathbf{V}(x)=0$ everywhere on $\mathcal{C}$.
 \end{theorem}

For the proofs of the next Theorem \ref{t1.4} and Theorem \ref{t1.5} ({\em Korns II inequality for domains ``with boundary condition"}) we refer to the same sources \cite{Ci00,Du10} mentioned above.
 %
 \begin{theorem}\label{t1.4} Let $1<p<\infty$, $\Omega\subset\mathbb{R}^n$ be a domain with the
Lip\-shitz boundary. Then the inequality
 \begin{eqnarray}\label{e1.18}
\big\|\mathbf{U}\big|\mathbb{W}^1_p(\Omega)\big\|\leq M\big\|{\rm\bf Def}(\mathbf{U})\big|\mathbb{L}_p(\Omega)\big\|
 \end{eqnarray}
holds with some constant $M>0$ or, equivalently, the equality
 \begin{eqnarray}\label{e1.19}
\big\|\mathbf{U}\big|\mathbb{W}^1_p(\Omega)\big\|_0:=\big\|{\rm\bf Def}(\mathbf{U})\big|\mathbb{L}_p(\Omega)\big\|
 \end{eqnarray}
defines an equivalent norm on the space $\widetilde{\mathbb{W}}{}^1_p(\Omega)$.
 \end{theorem}
 %
 \begin{theorem}\label{t1.5} Let $1<p<\infty$, $\mathcal{C}\subset\mathbb{R}^n$ be a Lip\-shitz
hypersurface with boundary. Then the inequality
 \begin{eqnarray}\label{e1.20}
\big\|\mathbf{V}\big|\mathbb{W}^1_p(\mathcal{C})\big\|\leq M\big\|{\rm\bf Def}_\mathcal{C}(\mathbf{V})\big|\mathbb{L}_p(\mathcal{C})\big\|
 \end{eqnarray}
holds with some constant $M>0$ or, equivalently, the equality
 \begin{eqnarray}\label{e1.21}
\big\|\mathbf{V}\big|\mathbb{W}^1_p(\mathcal{C})\big\|_0:=\big\|{\rm\bf Def}_\mathcal{C}(\mathbf{V})\big|\mathbb{L}_p(\mathcal{C})\big\|
 \end{eqnarray}
defines an equivalent norm on the space $\widetilde{\mathbb{W}}{}^1_p(\mathcal{C})$.
 \end{theorem}
 %
 \begin{remark}\label{r1.6}
A remarkable consequences of the foregoing theorems \ref{t1.2}-\ref{t1.5}
are the facts that the spaces $\mathbb{W}^1_p(\Omega)$ and $\widehat{\mathbb{W}}^1_p(\Omega)$ (as well as the spaces $\mathbb{W}^1_p(\mathcal{C})$ and $\widehat{\mathbb{W}}^1_p(\mathcal{C})$), where
 \begin{eqnarray*}
\widehat{\mathbb{W}}^1_p(\Omega):=\left\{\mathbf{U}=\big(U_1\ldots,U_n\big)^\top\,:\,
     U_j\,,\,\mathbf{D}_{jk}(\mathbf{U})\in\mathbb{L}_p(\Omega)\;\;\mbox{\rm for all}\;j,k=1,\ldots n\right\},\\[2mm]
\widehat{\mathbb{W}}^1_p(\mathcal{C}):=\left\{\mathbf{V}=\big(V_1\ldots,V_n\big)^\top\,:\,
     V_j\,,\,\mathfrak{D}_{jk}(\mathbf{V})\in\mathbb{L}_p(\mathcal{C})\;\;\mbox{\rm for all}\;j,k=1,\ldots n\right\}
 \end{eqnarray*}
are isomorphic (i.e. can be identified), although only $\displaystyle\frac{n(n+1)}2<n^2$ linear combinations of the $n^2$ derivatives $\partial_jU_k$ (of derivatives $\mathcal{D}_jU_k$, respectively), $j,k=1,\ldots n$ are involved in the definition of the equivalent norms in \eqref{e1.14} and \eqref{e1.14} (of the norms in \eqref{e1.19} and \eqref{e1.21}, respectively).
 \end{remark}

The next Lemma \ref{l1.7} is a slight generalization of \cite[Theorem 6.28.2]{Tr72} proved there for $p=2$.
 %
\begin{lemma}\label{l1.7} Let $\Omega$ be a bounded domain with the Lipschitz boundary ( (a
surface with the {\em uniform cone property}), $m=1,2,\ldots$, $1\leqslant p<\infty$ and let $F(\varphi)$ be a non-negative continuous functional on the Sobolev space $\mathbb{W}^m(\Omega)$:
\begin{itemize}
\item[i.]
$F\;:\;\mathbb{W}^m(\Omega)\to\mathbb{R}$ and $F(\lambda\varphi)=|\lambda| F(\varphi)$ for all complex $\lambda\in\mathbb{C}$ and all functions $\varphi\in\mathbb{H}^m(\Omega)$;
\item[ii.]
$0\leqslant F(\varphi)\leqslant C\big\|\varphi\big|\mathbb{W}^m_p(\Omega)\big\|$ for some constant $C>0$ and $F(P)\not=0$ for all polynomials of degree less than $m$.
\end{itemize}

Then the formula
 \begin{eqnarray}\label{e1.22}
\big\|\varphi\big|\mathbb{W}^m_p(\Omega)\big\|_F:=\big[\sum_{\alpha=m}
     \big\|\partial^\alpha\varphi\big|\mathbb{L}_p(\Omega)\big\|^p
     +F^p(\varphi)\big]^{1/p}
 \end{eqnarray}
defines an equivalent norm on the Sobolev space $\mathbb{W}^m(\Omega)$.

Lemma is valid if e replace $\Omega$ by a hypersurface $\mathcal{C}$ and partial derivatives $\partial^\alpha$-by G\'unters derivatives $\mathcal{D}^\alpha$.
 \end{lemma}
 \noindent
{\bf Proof:}
Let us note that $\big\|\varphi\big|\mathbb{W}^m_p(\Omega)\big\|_F$ in \eqref{e1.22} defines a norm on $\mathbb{W}^m_p(\Omega)$ indeed. Since other properties are trivial to check, we will only check that  $\big\|\varphi\big|\mathbb{W}^m_p(\Omega)\big\|_F=0$ implies $\varphi=0$. Then $F(\varphi)=0$ and all derivatives of order $m$ vanish: $\partial^\alpha\varphi=0$ for all $|\alpha|=m$. The latter means that the corresponding function is polynomial of order less than $m$, i.e., $\varphi(x)=\displaystyle\sum_{|\beta|<m}c_\beta x^\beta$. Since $F(\varphi)=0$, we get $\varphi\equiv0$ due to the property (ii).

Due to the condition (ii) holds the inequality
 \[
\|\psi\,\big|\,\mathbb{W}^m_p(\Omega)\|_F\leqslant\big[\sum_{\alpha=m}
     \big\|\partial^\alpha\varphi\big|\mathbb{L}_p(\Omega)\big\|^p
     +\|\psi\,\big|\,\mathbb{W}^m_p(\Omega)\|^p\big]^{1/p}\leqslant2^{1/p}
     \|\psi\,\big|\,\mathbb{W}^m_p(\Omega)\|.
 \]
Therefore the embedding of the spaces ${\mathbb{W}}^1_p(\Omega)\subset{\mathbb{W}}^1_{p,F}(\Omega)$, where
${\mathbb{W}}^1_{p,F}(\Omega)$ is the closure of $C^m(\Omega)$ with respect to the norm $\|\psi\,\big|\,\mathbb{W}^m_p(\Omega)\|_F$, is continuous.

If we apply the open mapping theorem of Banach (see \cite[Theorem 2.11, Corollary 2.12.b]{Ru73}, we conclude that the inverse inequality
 \[
\|\psi\,\big|\,\mathbb{W}^m_p(\Omega)\|\leqslant C\|\psi\,{\mathbb{W}}^m_{p,F}(\Omega)\,\|
     :=C\|\nabla\psi\,\big|\,\mathbb{W}^m_p(\Omega)\,\|_F
 \]
holds and accomplishes the proof. \hfill $\Box$

 %
 %
\section{Proofs of the basic inequalities}
 \label{sec2}
\setcounter{equation}{0}

\noindent
{\bf Proof of inequality \eqref{e0.1P}:} Let ${\mathbb{W}}^1_{p,\#}(\Omega)$ denote the subspace of ${\mathbb{W}}^1_{p,\#}(\Omega)$, consisiting of functions with mean value zero:
 \begin{eqnarray}\label{e2.1}
\varphi_\Omega:=\frac1{{\rm mes}\,\Omega}\int_\Omega\varphi(y)dy=0.
 \end{eqnarray}
The formula
\begin{eqnarray}\label{e2.2}
\|\varphi\,\big|\,{\mathbb{W}}^1_{p,\#}(\Omega)\,\|:=\|\nabla\varphi\,\big|\,\mathbb{L}_p(\Omega)\|
\end{eqnarray}
defines an equivalent norm in the space ${\mathbb{W}}^1_{p,\#}(\Omega)$. Since other properties are trivial to check, we only have to check that $\|\varphi\,\big|\,{\mathbb{W}}^1_{p,\#} (\Omega)\,\|=\|\nabla\varphi\,\big|\,\mathbb{L}_p(\Omega)\,\|=0$ implies $\varphi=0$. Indeed, the trivial norm implies that the gradient vanishes $\nabla\varphi=0$, which means that the corresponding function is constant $\varphi=C_0={\rm const}$; since the mean value is zero $\varphi_\Omega=C_0=0$ and $\varphi\equiv0$.

The inequality $\|\psi\,\big|\,{\mathbb{W}}^1_{p,\#}(\Omega)\,\|$ $\leqslant\|\psi\,\big|\, \mathbb{W}^1_p(\Omega)\|$, where
 \[
\|\psi\,\big|\,\mathbb{W}^1_p(\Omega)\|:=\left[\|\psi\,\big|\,\mathbb{L}_p(\Omega)\|^p
     +\|\nabla\psi\,\big|\,\mathbb{L}_p(\Omega)\|^2\right]^{1/p}
 \]
is the standard subspace norm on ${\mathbb{W}}^1_{p,\#}(\Omega)$ is trivial. Therefore the embedding
${\mathbb{W}}^1_{p,\#}(\Omega)\subset{\mathbb{W}}^1_p(\Omega)$ with the appropriate norms is continuous and proper, since constants belong to ${\mathbb{W}}^1_p(\Omega)$ but not to ${\mathbb{W}}^1_{p,\#}(\Omega)$.

If we apply the open mapping theorem of Banach (see \cite[Theorem 2.11, Corollary 2.12.b]{Ru73}, we conclude that the inverse inequality
 \[
\|\psi\,\big|\,\mathbb{W}^1_p(\Omega)\|\leqslant C_1\|\psi\,{\mathbb{W}}^1_{p,\#}(\Omega)\,\|
     =C_1\|\nabla\psi\,\big|\,\mathbb{L}_p(\Omega)\,\|
 \]
holds with some constant $C_1<\infty$ for all $\psi\in{\mathbb{W}}^1_{p,\#}(\Omega)$ (see \cite[Theorem 6.28.2]{Tr72} for a similar proof).

Since $\varphi_0:=\varphi-\varphi_\Omega\in{\mathbb{W}}^1_{p,\#}(\Omega)$, we have
 \[
\|\varphi-\varphi_\Omega\,\big|\,\mathbb{W}^1_p(\Omega)\|^p=\|\varphi-\varphi_\Omega\,\big|\,
     \mathbb{L}_p(\Omega)\,\|^p+\|\nabla\varphi\,\big|\,\mathbb{L}_p(\Omega)\,\|^p\leqslant C^p_1\|\nabla\varphi\,\big|\,\mathbb{L}_p(\Omega)\,\|^p.
 \]
The claimed inequality \eqref{e0.1P} follows with the constant $C:=(C_1^p-1)^{1/p}$.   \hfill $\Box$

\vskip3mm
\noindent
{\bf Proof of inequalities \eqref{e0.1F}, \eqref{e0.3}, \eqref{e0.14P} and \eqref{e0.14P2}:} Inequalities \eqref{e0.3} and \eqref{e0.14P2} are particular cases of \eqref{e0.14P}. Inequality \eqref{e0.14P} follows from \eqref{e1.22} if the functional $F$ is chosen as follows:
 \[
F(\varphi):=\left[\sum_{|\beta|<m}\left|\int_{\mathcal{M}_0}(\partial^\beta\varphi)^+(\cx)\,d\sigma
     \right|^p\right]^{1/p}.
 \]
The condition  $F(\varphi)\leqslant C\big\|\varphi\big|\mathbb{W}^m_p(\Omega)\big\|$ (see Lemma \ref{l1.7}.ii) holds due to the  Sobolev's continuous embeddings $\mathbb{W}^{m-k-1/p}_p(\mathcal{M}_0)\subset\mathbb{L}_p(\mathcal{M}_0)$, $k=0,1,\ldots,m-1$, and the trace theorem
 \[
\|(\partial^\beta\varphi)^+|\mathbb{W}^{m-|\beta|-1/p}_p(\mathcal{M}_0)\|\leqslant C_1\|\varphi|\mathbb{W}^m_p(\Omega)\|, \qquad |\beta|<m.
 \]
\vskip-7mm  \hfill $\Box$ \vskip4mm

\vskip3mm
\noindent
{\bf Proof of inequalities \eqref{e0.6}, \eqref{e0.6F}, \eqref{e0.7}, \eqref{e0.14P1} and \eqref{e0.14P3}:} Inequality \eqref{e0.6} is proved verbatim to \eqref{e0.1P}.

Inequalities \eqref{e0.6F}, \eqref{e0.7} and \eqref{e0.14P3} follow from \eqref{e0.14P1} (are particular cases).  Inequality \eqref{e0.14P1} follows from \eqref{e1.22} for surfaces if the functional $F$ is chosen as follows:
 \[
F(\varphi):=\left[\sum_{|\beta|<m}\left|\int_{\Gamma_0}(\mathcal{D}^\beta\varphi)^+(\cx)\,d\sigma
     \right|^p\right]^{1/p}.
 \]
The condition  $F(\varphi)\leqslant C\big\|\varphi\big|\mathbb{W}^m_p(\mathcal{C})\big\|$ (see Lemma \ref{l1.7}.ii) holds due to the  Sobolev's continuous embeddings $\mathbb{W}^{m-k-1/p}_p(\Gamma_0)\subset\mathbb{L}_p(\Gamma_0)$, $k=0,1,\ldots,m-1$, and the trace theorem
 \[
\|(\mathcal{D}^\beta\varphi)^+|\mathbb{W}^{m-|\beta|-1/p}_p(\Gamma_0)\|\leqslant C_1\|\varphi
     |\mathbb{W}^m_p(\mathcal{C})\|, \qquad |\beta|<m
 \]
\vskip-7mm  \hfill $\Box$\vskip4mm

\noindent
{\bf Proof of inequalities \eqref{e0.5} and  \eqref{e0.5a}:} Let $\Omega:=\mathcal{C}\times[a,b]$ and $\mathcal{M}_0:=\Gamma_0\times[a,b]$. To prove the inequality \eqref{e0.5} we proceed similarly: the
formula
 \[
\|\varphi\big|\mathbb{W}^1_p(\Omega)\|_0:=\left[\|\nabla_\mathcal{C}\varphi\big|\mathbb{L}_p(\Omega)\|^p
     +\displaystyle\int_{\mathcal{M}_0}\left|\varphi^+(\cx)\right|^pd\sigma\right]^{1/p}
 \]
defines a norm in the space $\mathbb{W}^1_p(\mathcal{C})$. Indeed, we have to check that $\|\varphi\big|\mathbb{W}^1_p(\Omega)\|_0=0$, which implies
\begin{eqnarray}\label{e2.2a}
\nabla_\mathcal{C}\varphi(\cx,t)=0 \quad \forall\,\cx\in\mathcal{C},\quad t\in[a,b],  \qquad
     \int_{\mathcal{M}_0}|\varphi^+(\tau,t)|^p\,d\sigma=0,
\end{eqnarray}
gives $\varphi=0$. But from the first equality in \eqref{e2.2a}, due to Lemma \ref{l1.1}, follows $\varphi(\cx,t)=\varphi(t)$ is independent of the surface variable.  But since $\varphi(\tau,t)=0$ on $\mathcal{M}_0$ (see the second equality in \eqref{e2.2a}), the function vanishes on the entire level surface $\mathcal{C}\times\{t\}$ for all $t\in[a,b]$. Then $\varphi=0$ in $\Omega$.

Due to Sobolev's continuous embedding $\mathbb{W}^{1-1/p}_p(\mathcal{M}_0)\subset\mathbb{L}_p(\mathcal{M}_0)$ and the trace theorem
\begin{eqnarray}\label{e2.2b}
\|\varphi^+|\mathbb{W}^{1-1/p}_p(\mathcal{M}_0)\|\leqslant C_1\|\varphi|\mathbb{W}^1_p(\Omega)\|
\end{eqnarray}
(see \cite{Tr95}), the initial norm in the space $\mathbb{W}^1_p(\Omega)$
 \[
\|\varphi\,\big|\,\mathbb{W}^1_p(\Omega)\|:=\left[\|\varphi\,
     \big|\,\mathbb{L}_p(\Omega)\|^p+\|\partial_t\varphi\,\big|\,\mathbb{L}_p(\Omega)\|^p
     +\|\nabla_\mathcal{C}\varphi\,\big|\,\mathbb{L}_p(\Omega)\|^2\right]^{1/p}
 \]
(see \cite{Du10,Du11}) estimates, obviously, the introduced norm
 \[
\|\nabla_\mathcal{C}\varphi\,\big|\,\mathbb{W}^1_p(\Omega)\|_0\leqslant C_2\|\varphi\,\big|\,\mathbb{W}^1_p(\Omega,\mathcal{M}_0)\|.
 \]
Then from open mapping theorem of Banach follows the inverse inequality
 \[
\|\varphi\,\big|\,\mathbb{L}_p(\Omega)\|\leqslant\|\varphi\,\big|\,\mathbb{W}^1_p(\Omega)\|\leqslant C
     \|\varphi\,\big|\,\mathbb{W}^1_p(\Omega)\|_0
 \]
and this accomplishes the proof of \eqref{e0.5}.

The inequality \eqref{e0.5a} is a direct consequence of \eqref{e0.5}.      \hfill $\Box$

\vskip3mm
\noindent
{\bf Proof of inequalities \eqref{e0.6}:} The proof is verbatim to the proof of inequality \eqref{e0.1P}, using the standard norm \eqref{e1.8} and the equivalent norm $\|\nabla_\mathcal{C}\varphi\big|\mathbb{L}_p(\mathcal{C})\|$ on the space ${\mathbb{W}}^1_{p,\#}(\mathcal{C})$. We also have to apply Lemma \ref{l1.1} to conclude that $\nabla_\mathcal{C}\varphi=0$ for $\varphi\in\mathbb{W}^1_p(\mathcal{C},\Gamma_0)$ implies $\varphi\equiv0$.  \hfill $\Box$

\vskip3mm
\noindent
{\bf Proof of inequalities \eqref{e0.6F}, \eqref{e0.7} and \eqref{e0.14P1}:} Inequality \eqref{e0.7} is a particular case of \eqref{e0.6F} (and of \eqref{e0.6}), while \eqref{e0.6F} is, in its turn, a particular case, $m=1$, of \eqref{e0.14P1}. Inequality \eqref{e0.14P1} follows from \eqref{e1.22} if the functional $F$ is chosen as follows:
 \[
F(\varphi):=\left[\sum_{|\beta|<m}\left|\int_{\mathcal{M}_0}(\partial^\beta\varphi)^+(\cx)\,d\sigma
     \right|^p\right]^{1/p}.
 \]
The condition  $F(\varphi)\leqslant C\big\|\varphi\big|\mathbb{W}^m_p(\Omega)\big\|$ (see Lemma \ref{l1.7}.ii) holds due to the  Sobolev's continuous embeddings $\mathbb{W}^{m-1/p}_p(\mathcal{M}_0)\subset\mathbb{L}_p(\mathcal{M}_0)$ and the trace theorem
 \[
\|\varphi^+|\mathbb{W}^{m-1/p}_p(\mathcal{M}_0)\|\leqslant C_1\|\varphi|\mathbb{W}^m_p(\Omega)\|.
 \]
(see \cite{Tr95}).    \hfill $\Box$

\vskip3mm
\noindent
{\bf Proof of inequalities \eqref{e0.14F} and \eqref{e0.14F1}:} These inequalities follow from \eqref{e1.22} if the functional $F$ is chosen as follows
 \[
F(\varphi):=\left[\int_{\mathcal{M}_0}\left|\varphi^+(\cx)\right|^p d\sigma
     \right]^{1/p}
 \]
for a domain $\Omega$ and
 \[
F(\varphi):=\left[\int_{\Gamma_0}\left|\varphi^+(\cx)\right|^p d\sigma
     \right]^{1/p}
 \]
 for a hypersurface $\mathcal{C}$ (see \eqref{e2.2b} for the justification of the condition  $F(\varphi)\leqslant C\big\|\varphi\big|\mathbb{W}^m_p\big\|$ in Lemma \ref{l1.7}.ii).       \hfill $\Box$

\vskip3mm
\noindent
{\bf Proof of inequality \eqref{e0.8} (and of similar ones):} For the space of smooth functions $C^1(\Omega)$ the proof is verbatim to the cases of the space $\mathbb{W}^1_p(\Omega)$.        \hfill $\Box$

\vskip3mm
\noindent
{\bf Proof of inequalities \eqref{e0.11} and  \eqref{e0.12}}: Based on the unique continuation property (see Theorem \ref{t1.2} and Theorem \ref{t1.3}), we prove easily that
 \begin{eqnarray*}
&\hskip-10mm \big\|\mathbf{U}\big|\mathbb{W}^1_p(\Omega)\big\|_1:=\left[\|{\rm\bf Def}\,\mathbf{U}\big|
    \mathbb{L}_p(\Omega)\|^p+\|\mathbf{U}^+\big|\mathbb{L}_p(\mathcal{M}_0)\|^p\right]^{1/p},\\[2mm]
&\hskip-10mm\|\mathbf{U}\big|\mathbb{W}^1_p(\mathcal{C})\|_1:=\left[\|{\rm\bf Def}_\mathcal{C}\mathbf{U}\big|
    \mathbb{L}_p(\mathcal{C})\|^p+\|\mathbf{U}^+\big|\mathbb{L}_p(\Gamma_0)\|^p\right]^{1/p}
 \end{eqnarray*}
define norms in the spaces $\mathbb{W}^1_p(\Omega)$ and in $\mathbb{W}^1_p(\mathcal{C})$, respectively. Then the obvious inequalities
 \begin{eqnarray*}
\|\mathbf{U}\big|\mathbb{W}^1_p(\Omega)\|_1\leqslant\big\|\mathbf{U}\big|\mathbb{W}^1_p(\Omega)\big\|_0\quad{\rm and}\quad
\|\mathbf{U}\big|\mathbb{W}^1_p(\mathcal{C})\|_1\leqslant\big\|\mathbf{U}\big|\mathbb{W}^1_p(\mathcal{C})\big\|_0
 \end{eqnarray*}
with the equivalent norms on the spaces $\mathbb{W}^1_p(\Omega)$ and in $\mathbb{W}^1_p(\mathcal{C})$ defined in \eqref{e1.14} and \eqref{e1.17} (see Theorem \ref{t1.4} and Theorem \ref{t1.5}) and the open mapping theorem of Banach ensure that the inverse inequalities
 \begin{eqnarray*}
\|\mathbf{U}\big|\mathbb{L}_p(\Omega)\|\leqslant C_1\|\mathbf{U}\big|\mathbb{W}^1_p(\Omega)\|_0\leqslant
     C\big\|\mathbf{U}\big|\mathbb{W}^1_p(\Omega)\big\|_1,\\[2mm]
\|\mathbf{U}\big|\mathbb{L}_p(\mathcal{C})\|\leqslant C_1\|\mathbf{U}\big|\mathbb{W}^1_p(\mathcal{C})\|_0\leqslant
     C\big\|\mathbf{U}\big|\mathbb{W}^1_p(\mathcal{C})\big\|_1.
 \end{eqnarray*}
hold and accomplish the proof.   \hfill $\Box$

\vskip3mm
\noindent
{\bf Proof of inequalities \eqref{e0.9} and  \eqref{e0.10}}: These inequalities are obvious consequences of \eqref{e0.11} and  \eqref{e0.12}.      \hfill $\Box$

\vskip3mm
\noindent
{\bf Proof of inequalities \eqref{e0.15} and  \eqref{e0.16}}: Inequality \eqref{e0.15} is proved verbatim to inequality \eqref{e0.5} by using, instead of Lemma \ref{l1.1}, the unique continuation property of Killing's vector fields, solutions to the equtions system ${\rm\bf Def}_\mathcal{C}\mathbf{U}=0$ (see Theorem \ref{t1.3}).

Inequality \eqref{e0.16} is an obvious consequence of \eqref{e0.15}.       \hfill $\Box$

 \baselineskip=12pt

\end{document}